\documentclass[12pt,a4paper]{amsart}
\usepackage[utf8]{inputenc}
\usepackage{amsmath, amsthm, amssymb, enumerate}

\allowdisplaybreaks

\usepackage{cleveref}
\crefname{section}{§}{§§}
\Crefname{section}{§}{§§}

\usepackage{xcolor}

\usepackage[utf8]{inputenc}
\usepackage{eepic,epic}
\usepackage{esint}
\usepackage{multirow}
\usepackage{array}
\usepackage[utf8]{inputenc}
\usepackage{graphicx}
\usepackage{caption}
\usepackage{booktabs}
\usepackage{siunitx}
\usepackage{makecell}
\usepackage{mathrsfs}
\usepackage{enumitem}

\def\cwedge{\bigcirc\kern-1.07em\wedge\ }

\newcommand \Tr{\operatorname{Tr}}

\newcommand \tR{\overline{R}}

\newtheorem{thm}{Theorem}[section]
\newtheorem{lem}[thm]{Lemma}
\newtheorem{prop}[thm]{Proposition}

\numberwithin{equation}{section}

\usepackage{float}

\theoremstyle{definition}

\newtheorem*{defn*}{Definition}

\newtheorem{ex}[thm]{Example}
\newtheorem*{rem*}{Remark}

\theoremstyle{plain}

\newtheorem*{thm*}{Theorem}
\newtheorem*{lem*}{Lemma}
\newtheorem*{theorema}{Theorem A}
\newtheorem*{theoremb}{Theorem B}
\newtheorem*{prop*}{Proposition}

\numberwithin{equation}{section}

\begin{document}
	
\title[Curvature identities for Einstein
manifolds of dimension 5 and 6]{Curvature identities for Einstein
manifolds \\of dimension 5 and 6}
   \author{Yunhee Euh, Jihun Kim, and
JeongHyeong Park}
\address{Department of Mathematics, Sungkyunkwan University, Suwon, 16419, Korea}
\email{prettyfish@skku.edu, jihunkim@skku.edu, parkj@skku.edu}

\thanks{This work was supported by the National Research Foundation of Korea(NRF) grant funded by the\\
\indent Korea government(MSIT) (NRF-2019R1A2C1083957).\\
\indent }

\maketitle
\begin{abstract}
Patterson discussed the curvature identities on Riemannian manifolds in \cite{Pa},
and a curvature identity for any 6-dimensional Riemannian manifold was independently
derived from the Chern-Gauss-Bonnet Theorem \cite{EPS17}. In this paper, we provide the explicit formulae of Patterson’s curvature identity that holds on 5-dimensional and 6-dimensional Einstein manifolds. We confirm that the curvature identities on the Einstein manifold from the previous work \cite{EPS17} are the same as the curvature identities deduced from Patterson’s result. We also provide examples that support the theorems.
\end{abstract}

\section{Introduction}


Let $M=(M,g)$ be an $m$-dimensional Riemannian manifold and $\nabla$ be the Levi-Civita connection of $g$. The curvature tensor $R$ on $M$ is defined by $R(X, Y)Z=[\nabla_X, \nabla_Y]Z-\nabla_{[X,Y]}Z$ for $X$, $Y$, $Z\in \mathfrak{X}(M)$, where $\mathfrak{X}(M)$ is the Lie algebra of all smooth vector fields on $M$. The Ricci tensor of $M$ is defined by $\rho(X,Y)=\Tr(Z\rightarrow R(Z,X)Y)$ and the scalar curvature of $M$ is obtained by $\tau = \Tr \rho$. Throughout the paper, we assume that the components of the tensor fields are with respect to a local orthonormal frame $\{e_i\}$ and we also adopt the Einstein convention on sum over repeated indices unless otherwise specified.

An $m$-dimensional Riemannian manifold $(M,g)$ is said to be {\it Einstein} if $\rho = \frac{\tau}{m}g$. The $m$-dimensional Einstein manifold $(M,g)$ is said to be {\it super-Einstein} if the following condition is satisfied on $M$
\begin{equation}\label{eq:supEin}
    {\tR}(X,Y)= \sum_{a, b, c =1}^{m} R(X,e_{a},e_{b},e_{c})R(Y,e_{a},e_{b},e_{c})=\dfrac{||R||^2}{m} g(X,Y)
\end{equation} with constant $||R||^2$ (see \cite{BV,EPS13,GW}). The condition \eqref{eq:supEin} has some geometric meanings: For a compact manifold, an Einstein metric is critical for the functional $\int_{M}||R||^2 dv_{g}$ restricted to $vol(M)=1$ if and only if $\tR(X,Y) = \frac{||R||^2}{m}g(X,Y)$ (see \cite[Corollary~4.72]{Bes1}). 
{{Boeckx and Vanhecke \cite{BV} showed that an Einstein manifold $M$ is super-Einstein if and only if the unit tangent sphere bundle $T_1 M$ equipped with the standard contact metric structure has constant scalar curvature (\cite[Proposition~3.6]{BV}).}}

Let $M_{q}(r,f)$ denote the mean-value of a real-valued function $f$ over a geodesic sphere $S(q;r)$ with center $q$ and radius $r$ in an $m$-dimensional Riemannian manifold $M$:
\begin{equation*}\label{avr}
M_{q}(r,f) := \frac{1}{\int_{S(q;r)} dv_{S(q;r)}}\, \int_{p\in S(q;r)} f(p) \  d v_{S(q;r)}.
\end{equation*} 
Gray and Willmore \cite{GW} showed the mean-value properties for an Einstein and super-Einstein manifold: They proved, from the expansions of $M_{q}(r,f)$, that the harmonic function $f$ near $q$ has the mean-value properties
\begin{equation*}
  M_{q}(r,f) = f(q) + O(r^6),\quad {\text{as}}\quad r\rightarrow 0
\end{equation*}
for an Einstein manifold (\cite[Theorem~1.1]{GW}) and
\begin{equation*}
  M_{q}(r,f) = f(q) + O(r^8),\quad {\text{as}}\quad r\rightarrow 0
\end{equation*}
for a super-Einstein manifold (\cite[Theorem~6.1]{GW}).

A Riemannian manifold {{$M$}} is said to be {\it $2$-stein} if there exist two functions $\mu_1, \mu_2$ on $M$ such that $\Tr R_X = \mu_1 ||X||^2$ and $ \Tr (R_X^2) = \mu_2 ||X||^4$, for all $p \in M$ and all $X \in T_pM$. Here, the Jacobi operator $R_X$ is defined by $R_X Y=R(Y,X)X$ for a tangent vectors $X$, $Y$ at a point $p\in M$. A unit vector field $V$ on $M$ is said to be a {\it harmonic vector field} if it is a critical point for the energy functional in
the set of all unit vector fields of $M$ \cite{Wood}. A contact
metric manifold whose characteristic vector field $\xi$ is a
harmonic vector field is called an {\it H-contact manifold}. Nikolayevsky and Park \cite{NP} showed that for a Riemannian manifold $M$, $T_1 M$ equipped with the standard contact metric structure is $H$-contact if and only if $M$ is 2-stein. Gilkey, Swann, and Vanhecke \cite{GSV} showed that a 4-dimensional manifold $M$ is 2-stein if and only if locally there is a choice of orientation of $M$ for which the metric is self-dual and Einstein (\cite[Theorem~2.6]{GSV}). From the definition of a 2-stein manifold, we can derive the super-Einstein conditions (see \cite[Chap.~6, \S E]{Bes2}). Thus, a 2-stein manifold is necessarily super-Einstein.

Euh, Park, and Sekigawa \cite{EPS13} dervied a curvature identity on any 4-dimensional manifold from the Chern-Gauss-Bonnet theorem. There are many applications of this identity (see \cite{EPS11,EPS013,GHMV}). On the other hand, Deszcz, Hotlo\'s, and Sent\"urk \cite{DHS} gave some curvature properties of 4-dimensional semi-Riemannian manifolds as an application of Patterson's curvature identity.

In this paper, we introduce the curvature identities on some $5$- and $6$-dimensional Riemannian manifolds such as Einstein, super-Einstein manifolds. Our main results are the following.

\smallskip

{{\begin{theorema}
Let $M=(M,g)$ be a 5-dimensional Riemannian manifold.
\begin{enumerate}[label={{\rm{(\alph*)}}}]
    \item\label{eq:5-einsteinid} If $M$ is Einstein, then the following curvature identity holds on $M${\rm :}
    \begin{equation*}
2\tau\tR_{ij}+4\check{R}_{ij}+4\hat{R}_{ij}-8\mathring{R}_{ij} = \Big(\frac{\tau}{5}||R||^2 + \frac{\tau^3}{25}\Big)g_{ij}.
\end{equation*}
    \item\label{eq:5-supEinid} If $M$ is super-Einstein, then the following curvature identity holds on $M${\rm :}
    {{\begin{equation*}
   4\mathring{R}_{ij} - 2\hat{R}_{ij} =\Big(\frac{9}{50}\tau||R||^2 -\frac{\tau^3}{50}\Big)g_{ij}.
\end{equation*}}}
\end{enumerate}
\end{theorema}

\smallskip

\begin{theoremb}
Let $M=(M,g)$ be a 6-dimensional Riemannian manifold.
\begin{enumerate}[label={{\rm{(\alph*)}}}]
    \item\label{eq:6-einstein} If $M$ is Einstein, then the following curvature identity holds on $M${\rm :}
    {{\begin{equation*}
      4\tau {{\tR_{ij}}}+12\check{R}_{ij}+12\hat{R}_{ij}-24\mathring{R}_{ij}=(\tau ||R||^2 -4\mathring{R}+2\hat{R})g_{ij}.
   \end{equation*}}}
   \item\label{eq:6-dim2stein} If $M$ is super-Einstein, then the following curvature identity holds on $M${\rm :}
   \begin{equation*}
				2\mathring{R}_{ij}-\hat{R}_{ij}=\frac{1}{6}(2\mathring{R}-\hat{R})g_{ij}.
	\end{equation*}
\end{enumerate}
Here, we set
\begin{gather*}
    \check{R}_{ij}=R_{iuvj}R_{abcu}R_{abcv},\quad \hat{R}_{ij}=R_{ibcd}R_{jbuv}R_{cduv},\quad \mathring{R}_{ij}=R_{ibcd}R_{jucv}R_{budv},\\
    \hat{R}= R_{abcd}R_{abuv}R_{cduv},\quad \mathring{R}=R_{abcd}R_{aucv}R_{budv}.
\end{gather*}
\end{theoremb}}}

We derive the curvature identities from Patterson's curvature identity. In Section 2, we recall the previous results for a 6-dimensional Riemannian manifold \cite{EPS17} and introduce Patterson's curvature identity. In Section 3 and Section 4, we prove Theorem A and Theorem B, respectively. We also give some examples supporting the theorems.
In the appendix, we attach the detailed computation for the proof of Lemma \ref{thm:6-einstien(0,6)-tensor}.

\section{{{Curvature identities}}}

On a 6-dimensional compact oriented Riemannian manifold, Chern-Gauss-Bonnet theorem states that Euler characteristic $\chi(M)$ of $M$ is given by the following integral formula.
   \begin{prop}[\cite{Sa}]
   Let $M=(M,g)$ be a $6$-dimensional compact oriented Riemannian manifold. Then
the Euler characteristic $\chi(M)$ of $M$ is given by
      \begin{equation*}\label{eq:id6}
         \begin{aligned}
            \chi(M)
            =\frac{1}{384\pi^3}\int_M\{&\tau^3-12\tau||\rho||^2+3\tau||R||^2+16\rho_{ab}\rho_{ac}\rho_{bc}\\
            &-24\rho_{ab}\rho_{cd}R_{acbd}
            -24\rho_{uv}R_{abcu}R_{abcv}\\
             &+8R_{abcd}R_{aucv}R_{bvdu}-2R_{abcd}R_{abuv}R_{cduv}\}dV_g.
         \end{aligned}
      \end{equation*}
 \end{prop}

In \cite{EPS17}, the authors consider the one-parameter deformation $g(t)$ of $g$, using the fact that
Euler characteristic is a topological invariant for the deformation, they got the universal
curvature identity which holds on any 6-dimensional Riemannian manifold. In particular, they have the following.

\begin{thm}[\cite{EPS17}]\label{th:6-eins}
    Let $M=(M,g)$ be a 6-dimensional Einstein manifold. Then the following identity holds on $M${\rm :}
    \begin{equation*}
      (-\tau ||R||^2+4\mathring{R}-2\hat{R})g_{ij}+12\check{R}_{ij}+12\hat{R}_{ij}-24\mathring{R}_{ij}+4\tau {{\tR_{ij}}}=0.
   \end{equation*}
\end{thm}
\smallskip

{{For an integer $N \ge 1$, the {\textit{generalized Kronecker delta}} is given by
\begin{equation*}
	    \delta^{j_1j_2\cdots j_N}_{i_1i_2\cdots i_N}=
	    \begin{vmatrix}
	    \delta^{j_1}_{i_1}&\delta^{j_2}_{i_1}&\cdots&\delta^{j_N}_{i_1}\\
	    \delta^{j_1}_{i_2}&\delta^{j_2}_{i_2}&\cdots &\delta^{j_N}_{i_2}\\
	    \vdots&\vdots&\vdots&\vdots\\
	    \delta^{j_1}_{i_N}&\delta^{j_2}_{i_N}&\cdots&\delta^{j_N}_{i_N}
	    \end{vmatrix}.
\end{equation*}
By using the skew-symmetric properties of the generalized Kronecker delta, Patterson [14] obtained the curvature identity for an $m$-dimensional Riemannian manifold as follows:}}
\begin{equation}\label{eq:Pa_n}
    \delta^{jj_1j_2\cdots j_m}_{ii_1i_2\cdots i_m}{R^{i_1i_2}}_{j_1j_2} \cdots {R^{i_{2r-1}i_{2r}}}_{j_{2r-1}j_{2r}}=0,
\end{equation}
where $r$ is any integer such that $1\leq r\leq \frac{m}{2}$ ($m$ even) or $1\leq r\leq \frac{m-1}{2}$ ($m$ odd).

We note that the identity~\eqref{eq:Pa_n} holds with respect to the Weyl curvature tensor $W$ even if we replace $R$ by $W$ in \eqref{eq:Pa_n} (\cite[Section~8]{Pa}).

\section{Proof of Theorem A}
Let $(M,g)$ be a 5-dimensional Riemannian manifold. Making use of \eqref{eq:Pa_n}, we have the curvature identity with respect to the Weyl curvature tensor $W$ for the case when $r=2$ as follows:
\begin{equation}\label{eq:5-Weyl}
    \begin{aligned}
    &||W||^2(g_{ik}g_{jl}-g_{il}g_{jk})-4\Big\{W_{abcj}W_{abcl}g_{ik}-W_{abcj}W_{abck}g_{il}+W_{abci}W_{abck}g_{jl}\\
    &-W_{abci}W_{abcl}g_{jk}-2W_{iabl}W_{kabj}+2W_{iabk}W_{labj}-W_{abij}W_{abkl}\Big\}=0.
    \end{aligned}
\end{equation}

Since the Weyl tensor $W$ on a 5-dimensional Riemannian manifold is given by
\begin{equation*}
W_{pqrs}=R_{pqrs}-\frac{1}{3}(\rho_{ps}g_{qr}+\rho_{qr}g_{ps}-\rho_{pr}g_{qs}-\rho_{qs}g_{pr})+\frac{\tau}{12}(g_{ps}g_{qr}-g_{pr}g_{qs}),
\end{equation*}
we obtain the explicit formula of the curvature identity on $5$-dimensional Einstein manifolds.
\begin{lem}
Let $M=(M,g)$ be a 5-dimensional Einstein manifold. The following curvature identity holds on $M${\rm :}
\begin{equation}\label{eq:5-einstein}
    \begin{aligned}
    &\Big(||R||^2 + \frac{\tau^2}{5} \Big)(g_{ik}g_{jl} - g_{il}g_{jk})-4(\tR_{ik}g_{jl} + \tR_{jl}g_{ik} -\tR_{il}g_{jk} -\tR_{jk}g_{il})\\
    &~+8(R_{iabl}R_{kabj}-R_{iabk}R_{labj}) +4R_{abij}R_{abkl} + \frac{12}{5}\tau R_{ijkl} = 0.
    \end{aligned}
\end{equation}
\begin{proof}
Since $\rho_{ps}=\frac{\tau}{5}g_{ps}$, we obtain $W_{pqrs}=R_{pqrs}-\frac{\tau}{20}(g_{ps}g_{qr}-g_{pr}g_{ps})$. Now we apply this $W$ to the identity~\eqref{eq:5-Weyl}. Then we have
\begin{gather*}
    ||W||^2=||R||^2-\frac{\tau}{10},\quad W_{abcj}W_{abcl}=R_{abcj}R_{abcl}-\frac{\tau^2}{50}g_{jl},\\
    W_{iabj}W_{kabj}=R_{iabl}R_{kabj}-\frac{\tau}{10}R_{ijlk}-\frac{\tau^2}{80}g_{il}g_{jk}+\frac{\tau^2}{400}g_{ik}g_{jl},\\
     W_{abij}W_{abkl}=R_{abij}R_{abkl}+\frac{\tau}{5}R_{ijkl}+\frac{\tau^2}{200}\big(g_{ik}g_{jl}-g_{il}g_{jk}\big).
\end{gather*}
Similarly, we can obtain remaining terms. Then, by rearranging all terms, we complete the proof.
\end{proof}
\end{lem}

{{We note that there are useful formulae as follows:
\begin{equation}\label{eq:re}
R_{ibcd}R_{jbuv}R_{cudv}=\frac{1}{2}\hat{R}_{ij},\quad R_{ibcd}R_{jubv}R_{cudv}=\frac{1}{4}\hat{R}_{ij},\quad R_{ibcd}R_{jucv}R_{bvdu}=\mathring{R}_{ij}-\frac{1}{4}\hat{R}_{ij}.
\end{equation}}}

Now, we transvect each term of ~\eqref{eq:5-einstein} with $R_{pjkl}$ and use $\rho_{ip}=\frac{\tau}{5}g_{ip}$. Then, we obtain
\begin{equation*}
    \begin{aligned}
    &\Big(||R||^2 + \frac{\tau^2}{5} \Big)(g_{ik}g_{jl} - g_{il}g_{jk})R_{pjkl} = -\frac{2}{5}\tau\Big(||R||^2 + \frac{\tau^2}{5}\Big)g_{ip},\\
    -&4(\tR_{ik}g_{jl} + \tR_{jl}g_{ik} -\tR_{il}g_{jk} -\tR_{jk}g_{il})R_{pjkl} = -4\Big(-\frac{2}{5}\tau\tR_{ip} - 2\check{R}_{ip}\Big),\\
    &~8(R_{iabl}R_{kabj}-R_{iabk}R_{labj})R_{pjkl} = -16\mathring{R}_{ip} + 4\hat{R}_{ip},\\
    &~4R_{abij}R_{abkl}R_{pjkl} = 4\hat{R}_{ip},\\
    &~\frac{12}{5}\tau R_{ijkl}R_{pjkl} = \frac{12}{5}\tau\tR_{ip}.
    \end{aligned}
\end{equation*}
For the third equation above, we use~\eqref{eq:re}. Hence, we complete the proof of Theorem~A-\ref{eq:5-einsteinid}.

\smallskip

{{Now we prove Theorem~A-\ref{eq:5-supEinid}. Applying the condition~\eqref{eq:supEin} to~\eqref{eq:5-einstein}, we obtain}}
\begin{equation}\label{eq:Pa5}
	        R_{ijab}R_{abkl}+2R_{iabl}R_{kabj}-2R_{iabk}R_{labj}+\frac{3}{5}\tau R_{ijkl}=\Big(\frac{3}{20}||R||^2-\frac{1}{20}\tau^2\Big) (g_{ik}g_{jl}-g_{il}g_{jk}).
	        \end{equation}
Transvecting each term of \eqref{eq:Pa5} with $R_{pjkl}$ and using~\eqref{eq:supEin}, we obtain
\begin{equation*}\label{eq:5-supEin}
\begin{aligned}
     &R_{ijab}R_{abkl}R_{pjkl}=\hat{R}_{ip},\\
    &2R_{iabl}R_{kabj}R_{pjkl}= -2\Big(\mathring{R}_{ip} - \frac{1}{4}\hat{R}_{ip}\Big),\\
    -&2R_{iabk}R_{labj}R_{pjkl}=-2\Big(\mathring{R}_{ip} - \frac{1}{4}\hat{R}_{ip}\Big),\\
    &\frac{3}{5}\tau R_{ijkl}R_{pjkl}=\frac{3}{25}\tau||R||^2g_{ip},\\
    &(g_{ik}g_{jl}-g_{il}g_{jk}) R_{pjkl}=-\frac{2}{5}\tau g_{ip}.
\end{aligned}
\end{equation*}
For the second and third equations, we use \eqref{eq:re}.
Therefore, we complete the proof of Theorem~A-\ref{eq:5-supEinid}.

\smallskip

For a 5-dimensional 2-stein manifold, there is the following orthonormal basis introduced by Nikolayevsky \cite{Nik}.
	
	\begin{prop}\big(\cite[Proposition~4]{Nik}\big)
	Let $M=(M,g)$ be a 2-stein manifold of dimension 5. Then there exists an orthonomal basis $\{e_i\}$ at each point $p\in M$ such that
	\begin{equation*}\label{eq:basis5}
	    \begin{aligned}
	        R_{1212}&=R_{1313}=R_{2323}=R_{2424}=R_{3434}=\alpha-\beta,~
	        R_{1414}=\alpha-4\beta,\\ R_{1515}&=R_{4545}=\alpha,~ R_{2525}=R_{3535}=\alpha-3\beta,\\
	        R_{1234}&=\beta,~R_{1235}=\sqrt{3}\beta,~R_{1324}=-\beta,~R_{1325}=\sqrt{3}\beta,\\
	        R_{1423}&=-2\beta, ~R_{2425}=\sqrt{3}\beta,~R_{3435}=-\sqrt{3}\beta,	
	        \end{aligned}
	\end{equation*}
	and all the other components of $R$ vanish.
	\end{prop}
	By using {{Nikolayevsky's}} basis, we show the following theorem.
	\begin{thm}\label{thm:3.4}
	Let $M=(M,g)$ be a 5-dimensional 2-stein manifold. The identity \eqref{eq:Pa5} holds on $M$.
	\begin{proof}
Each term of the left hand side of \eqref{eq:Pa5} in the case of $i=1$, $j=2$, $k=3$, $l=4$  is as follows:
\begin{equation*}
    \begin{aligned}
      R_{12ab}R_{ab34}&=2(2 \alpha - 5 \beta)\beta,\\
      2R_{1ab4}R_{3ab2}& = 2(2\alpha-5\beta)\beta,\\
      -2R_{1ab3}R_{4ab2}& = 2(2\alpha+\beta)\beta,\\
    \frac{3}{5}\tau R_{1234}& = -6(2\alpha-3\beta)\beta.\\
    \end{aligned}
\end{equation*}
Since $R_{12ab}R_{ab34}+2R_{1ab4}R_{3ab2}-2R_{1ab3}R_{4ab2}+\frac{3}{5}\tau R_{1234}=0$ and
$(g_{13}g_{24}-g_{14}g_{23})=0$ in \eqref{eq:Pa5}, the equation \eqref{eq:Pa5} holds for $i=1$, $j=2$, $k=3$, $l=4$. For other choice of $i, j, k, l$,
similar processes can be done for the proof of Theorem \ref{thm:3.4}.
\end{proof}
	\end{thm}

Now we give the examples for Theorem A.

\begin{ex}
Let $M$ be a Riemannian product manifold of a 3-dimensional Riemannian manifold $M_{1}$ of constant sectional curvature $k$ and a surface $M_{2}$ of constant sectional curvature $2k$ $(k\neq 0)$. Then $M$ is not super-Einstein but Einstein manifold. Let $\{e_{i}\}$, $i=1,\ldots,5$ be an orthonormal basis of $M$, where $\{e_{1},e_{2},e_{3}\}$ and $\{e_{4},e_{5}\}$ are bases for $M_1 $ and $M_2 $, respectively. Then we have
\begin{equation}\label{eq:ex_5EinnotsupEin}
    R_{1221}=R_{1331}=R_{2332}=k,\quad R_{4554}=2k,
\end{equation}
and all the other components of $R$ vanish. From \eqref{eq:ex_5EinnotsupEin} we have $\tau = 10k$, $||R||^2 = 28k^2$, and
\begin{gather*}
    \tR_{ij}=
    \begin{cases}
    4k^2 & {\text{if}}\;\, i=j=1,2,3\\
    8k^2 & {\text{if}}\;\, i=j=4,5\\
    0 & {\text{otherwise}}
    \end{cases},\qquad
    \check{R}_{ij}=
    \begin{cases}
    8k^3 & {\text{if}}\;\, i=j=1,2,3\\
    16k^3 & {\text{if}}\;\, i=j=4,5\\
    0 & {\text{otherwise}}
    \end{cases},\\
    \hat{R}_{ij}=
    \begin{cases}
    -8k^3  & {\text{if}}\;\, i=j=1,2,3\\
    -32k^3 & {\text{if}}\;\, i=j=4,5\\
    0 & {\text{otherwise}}
    \end{cases},\qquad
    \mathring{R}_{ij}=
    \begin{cases}
    -2k^3 &  {\text{if}}\;\, i=j=1,2,3\\
    0 & {\text{otherwise}}
    \end{cases}.
\end{gather*}
Therefore, we find that the curvature identity of Theorem ~A-\ref{eq:5-supEinid} does not hold, but that of Theorem~A-\ref{eq:5-einsteinid} holds on $M$.
\end{ex}

\begin{ex}
{{Let}} $M=SL(3)/SO(3)$. Then $M$ is a 5-dimensional 2-stein manifold.
 The inner product and the curvature tensor are given by
\begin{equation*}
    <X,Y>=\Tr{XY},\quad R(X,Y)Z= - [[X,Y],Z]
\end{equation*}
for $X$, $Y$, $Z\in \mathfrak{X}(M)$. Then, with orthonormal basis of $SL(3)/SO(3)$
\begin{equation*}
\begin{gathered}
X_1=\frac{1}{\sqrt{2}}\begin{pmatrix}
0 & 1 & 0 \\
1& 0 & 0 \\
0 & 0 & 0
\end{pmatrix}, \quad
X_2=\frac{1}{\sqrt{2}}\begin{pmatrix}
0&0&1\\
0&0&0\\
1&0&0
\end{pmatrix},\quad
X_3=\frac{1}{\sqrt{2}}\begin{pmatrix}
0&0&0\\
0&0&1\\
0&1&0
\end{pmatrix}\\
X_4=\frac{1}{\sqrt{2}}\begin{pmatrix}
-1&0&0\\
0&1&0\\
0&0&0
\end{pmatrix},\quad
X_5=\frac{1}{\sqrt{6}}\begin{pmatrix}
1&0&0\\
0&1&0\\
0&0&-2
\end{pmatrix},
\end{gathered}
\end{equation*}
we have the components of the curvature tensor:
\begin{equation}\label{eq:ex_5-2stein}
    \begin{aligned}
    &~R_{1221}=R_{1331}=R_{2332}=R_{2442}=R_{3443}= -\frac{1}{2},~
	        R_{1441}= - 2,\\
	        &~R_{1551}=R_{4554}=0,~R_{2552}=R_{3553}=-\frac{3}{2},\\
	        &~R_{1234}=-\frac{1}{2},~
	        R_{1235}=-\frac{\sqrt{3}}{2},~R_{1324}=\frac{1}{2},~R_{1325}=-\frac{\sqrt{3}}{2},\\
	        &~R_{1423}=1,~R_{2425}=-\frac{\sqrt{3}}{2},~R_{3435}=\frac{\sqrt{3}}{2},
    \end{aligned}
    \end{equation}
and all the other components vanish.
{{From \eqref{eq:ex_5-2stein}, we compute the Ricci tensor $\rho_{ij}$. For $i\neq j$, $\rho_{ij}=0$ and
$\rho_{11}=\sum\limits_{a=1}^{5}R_{1aa1}=-3$. Similarly, $\rho_{22}= \rho_{33}= \rho_{44}= \rho_{55}= - 3$}} and so we get the scalar curvature $\tau$ is $-15$. Making use of \eqref{eq:ex_5-2stein}, ${{||R||^2=75}}$, $\hat{R}_{ij}=75\delta_{ij}$, and $\mathring{R}_{ij}=\frac{15}{4}\delta_{ij}$, where $\delta_{ij}$ denotes the Kronecker delta. Therefore, we find that the curvature identity of Theorem~A-\ref{eq:5-supEinid} holds on $M$.
\end{ex}

\section{Proof of Theorem B}

We obtain the explicit formula of the curvature identity on 6-dimensional Einstein manifolds for the case when $m=6$ and $r=2$ of the identity~\eqref{eq:Pa_n}. To prove the Theorem B, we need Lemma \ref{thm:6-einstien(0,6)-tensor}, which is proved in the Appendix.
\begin{lem}\label{thm:6-einstien(0,6)-tensor}
	Let $M=(M,g)$ be a $6$-dimensional Einstein manifold. The following curvature identity holds on $M${\rm :}
	\begin{equation}\label{eq:6-einstein(0,6)-tensor}
		\begin{aligned}
			&\Big\{\frac{1}{8}\Big(||R||^2 +\frac{\tau^2}{3}\Big)\Big(g_{ij}g_{hk}g_{lm}-g_{ij}g_{hm}g_{lk}-g_{ik}g_{hj}g_{lm}+g_{ik}g_{hm}g_{lj}+g_{im}g_{hj}g_{lk}-g_{im}g_{hk}g_{lj}\Big)\Big\}\\
			&+ \Big\{-\frac{1}{2}\Big(\tR_{ij}(g_{hk}g_{lm}-g_{hm}g_{lk})-\tR_{ik}(g_{hj}g_{lm}-g_{hm}g_{lj})+\tR_{im}(g_{hj}g_{lk}-g_{hk}g_{lj})\\
			&\qquad-\tR_{hj}(g_{ik}g_{lm}-g_{im}g_{lk})+\tR_{hk}(g_{ij}g_{lm}-g_{im}g_{lj})-\tR_{hm}(g_{ij}g_{lk}-g_{ik}g_{lj})\\
			&\qquad+\tR_{lj}(g_{ik}g_{hm}-g_{im}g_{hk})-\tR_{lk}(g_{ij}g_{hm}-g_{im}g_{hj})+\tR_{lm}(g_{ij}g_{hk}-g_{ik}g_{hj})\Big)\Big\}\\
			&+\Big\{\Big(-T_{ijkh} + T_{ikjh} +\frac{1}{2}S_{ihjk}+\frac{\tau}{3} R_{ihjk}\Big)g_{lm}-\Big(-T_{ijkl} +T_{ikjl}+\frac{1}{2}S_{iljk}+\frac{\tau}{3} R_{iljk}\Big)g_{hm}\\
			&-\Big(-T_{ijmh}+T_{imjh}+\frac{1}{2}S_{ihjm}+\frac{\tau}{3} R_{ihjm}\Big)g_{lk}+\Big(-T_{ijml} +T_{imjl}+\frac{1}{2}S_{iljm}+\frac{\tau}{3} R_{iljm}\Big)g_{hk}\\
			&+\Big(-T_{ikmh} +T_{imkh}+\frac{1}{2}S_{ihkm}+\frac{\tau}{3} R_{ihkm}\Big)g_{lj}-\Big(-T_{ikml} +T_{imkl}+\frac{1}{2}S_{ilkm}+\frac{\tau}{3}R_{ilkm}\Big)g_{hj}\\
			&-\Big(-T_{hjml} +T_{hmjl}+\frac{1}{2}S_{hljm}+\frac{\tau}{3} R_{hljm}\Big)g_{ik}+\Big(-T_{hjkl} +T_{hkjl}+\frac{1}{2}S_{hljk}+\frac{\tau}{3} R_{hljk}\Big)g_{im}\\
			&+\Big(-T_{hkml} +T_{hmkl}+\frac{1}{2}S_{hlkm}+\frac{\tau}{3} R_{hlkm}\Big)g_{ij}\Big\}\\
			&+A_{hjkmil}-A_{ljkmih}-A_{ijkmhl}+A_{ijmkhl}-A_{hjmkil}+A_{ljmkih}-A_{ikmjhl}+A_{hkmjil}-A_{lkmjih}\\
			&=0,
		\end{aligned}
	\end{equation}
	where $T_{pqrs} = R_{pabq}R_{rabs}$, $S_{pqrs} = R_{abpq}R_{abrs}$, and $A_{pqrstu} = R_{apqr}R_{astu}$.
\begin{proof}
The proof is similar to that of Lemma 3.1 for the case when $m=6$ and $r=2$. Here, the Weyl tensor W is given by
\begin{equation*}
    W_{abcd} = R_{abcd} - \frac{1}{4}(\rho_{ad}g_{bc}+\rho_{bc}g_{ad} - \rho_{ac}g_{bd} -\rho_{bd}g_{ac}) + \frac{\tau}{20}(g_{ad}g_{bc} - g_{ac}g_{bd}).
\end{equation*}
Next we expand \eqref{eq:Pa_n} for $W$ and rearrange them. Using the Einstein condition $\rho_{ij} = \frac{\tau}{6}g_{ij}$, we obtain 34 terms. Here we deal with a few terms among them. For more details, we refer to appendix. The numbers shown below correspond to those in appendix.
\begin{enumerate}[label={(\arabic*)},leftmargin=0cm,align=left]
\item[(1)] $\big(||R||^2 - 4||\rho||^2 + \tau^2 \big)g_{ij}g_{hk}g_{lm}=\left(||R||^2 +\dfrac{\tau^2}{3} \right)g_{ij}g_{hk}g_{lm}$
\item[(7)]
$\big(-4R_{iabc}R_{jabc} + 8\rho_{ia}\rho_{ja} +8R_{iabj}\rho_{ab}-4\tau\rho_{ij}\big)g_{hk}g_{lm}=-4\tR_{ij}g_{hk}g_{lm} -\dfrac{2}{9}\tau^2 g_{ij}g_{hk}g_{lm}$
\item[(25)]
\begin{align*}
        &\big(-8R_{iabj}R_{kabh} +8R_{iabk}R_{jabh} +4R_{abih}R_{abjk} -8R_{aijk}\rho_{ah} +8R_{akih}\rho_{aj} -8R_{ajih}\rho_{ak}\\
        &\quad+8R_{ahjk}\rho_{ai} +8\rho_{ij}\rho_{hk} -8\rho_{ik}\rho_{hj} -4\tau R_{ihjk}\big)g_{lm}\\
        &=\Big(-8T_{ijkh} +8T_{ikjh} +4S_{ihjk} +\dfrac{4}{3}\tau R_{ihjk} +\dfrac{2}{9}\tau^2 g_{ij}g_{hk} - \dfrac{2}{9}\tau^2 g_{ik}g_{hj} \Big)g_{lm}
\end{align*}
\item[(34)]
\begin{align*}
        &8R_{ahjk}R_{amil}-8R_{aljk}R_{amih}-8R_{aijk}R_{amhl}+8R_{aijm}R_{akhl}-8R_{ahjm}R_{akil}+8R_{aljm}R_{akih}\\
        &\;-8R_{aikm}R_{ajhl}+8R_{ahkm}R_{ajil}-8R_{alkm}R_{ajih} -8R_{ilkm}\rho_{hj}+8R_{iljm}\rho_{hk}+8R_{hljk}\rho_{im}\\
        &\;+8R_{ihjk}\rho_{lm}+8R_{ihkm}\rho_{lj}-8R_{ihjm}\rho_{lk}-8R_{iljk}\rho_{hm}+8R_{hlkm}\rho_{ij}-8R_{hljm}\rho_{ik}\\
        &= 8(A_{hjkmil}-A_{ljkmih}-A_{ijkmhl}+A_{ijmkhl}-A_{hjmkil}+A_{ljmkih}-A_{ikmjhl}+A_{hkmjil}\\
        &\quad-A_{lkmjih}){{+\frac{4}{3}\tau}} \Big({{-}}R_{ilkm}g_{hj}+ R_{iljm}g_{hk}+ R_{hljk}g_{im}+ R_{ihjk}g_{lm}+ R_{ihkm}g_{lj}- R_{ihjm}g_{lk}\\
        &\quad- R_{iljk}g_{hm}+R_{hlkm}g_{ij}- R_{hljm}g_{ik}\Big).
\end{align*}
\end{enumerate}
By rearranging all terms, we complete the proof of Lemma~\ref{thm:6-einstien(0,6)-tensor}.
\end{proof}
\end{lem}

Now, we transvect {{each term}} of \eqref{eq:6-einstein(0,6)-tensor} with $R_{ihjk}$.
For the first term of \eqref{eq:6-einstein(0,6)-tensor}, we have
\begin{equation*}
\begin{aligned}
&\frac{1}{8}\Big(||R||^2 +\frac{\tau^2}{3}\Big)\Big(g_{ij}g_{hk}g_{lm}-g_{ij}g_{hm}g_{lk}-g_{ik}g_{hj}g_{lm}+g_{ik}g_{hm}g_{lj}+g_{im}g_{hj}g_{lk}-g_{im}g_{hk}g_{lj}\Big) R_{ihjk}\\
&=-\frac{1}{6}\Big(\tau||R||^2 +\frac{\tau^3}{3}\Big) g_{lm}.	
\end{aligned}
\end{equation*}
\noindent For the second term, we have
\begin{equation*}
-\frac{1}{2}\tR_{ij}(g_{hk}g_{lm}-g_{hm}g_{lk})R_{ihjk}=\frac{\tau}{12}||R||^2g_{lm}-\frac{1}{2}\check{R}_{lm},
\end{equation*}
\begin{equation*}
\frac{1}{2}\tR_{ik}(g_{hj}g_{lm}-g_{hm}g_{lj}) R_{ihjk}=\frac{\tau}{12}||R||^2g_{lm}-\frac{1}{2}\check{R}_{lm},\;    -\frac{1}{2}\tR_{im}(g_{hj}g_{lk}-g_{hk}g_{lj}) R_{ihjk}=-\frac{\tau}{6}\tR_{lm},
\end{equation*}
\begin{equation*}
\frac{1}{2}\tR_{hj}(g_{ik}g_{lm}-g_{im}g_{lk})R_{ihjk}=\frac{\tau}{12}||R||^2g_{lm}-\frac{1}{2}\check{R}_{lm},
\end{equation*}
\begin{equation*}
-\frac{1}{2}\tR_{hk}(g_{ij}g_{lm}-g_{im}g_{lj})R_{ihjk}=\frac{\tau}{12}||R||^2g_{lm}-\frac{1}{2}\check{R}_{lm},
\end{equation*}
\begin{equation*}
\frac{1}{2}\tR_{hm}(g_{ij}g_{lk}-g_{ik}g_{lj})R_{ihjk}=-\frac{\tau}{6}\tR_{lm}, \;-\frac{1}{2}\tR_{lj}(g_{ik}g_{hm}-g_{im}g_{hk}) R_{ihjk}=-\frac{\tau}{6}\tR_{lm},
\end{equation*}
\begin{equation*}
\frac{1}{2}\tR_{lk}(g_{ij}g_{hm}-g_{im}g_{hj}) R_{ihjk}=-\frac{\tau}{6}\tR_{lm},\; -\frac{1}{2}\tR_{lm}(g_{ij}g_{hk}-g_{ik}g_{hj})R_{ihjk}=\tau\tR_{lm}.
\end{equation*}
Thus, the second term of \eqref{eq:6-einstein(0,6)-tensor} transvecting with $R_{ihjk}$  becomes $\dfrac{\tau}{3}||R||^2 g_{lm}-2\check{R}_{lm}+\dfrac{\tau}{3}\tR_{lm}$.\\
For the third term, we have
\begin{equation*}
\begin{aligned}
&\Big(-T_{ijkh} + T_{ikjh} +\frac{1}{2}S_{ihjk}+\frac{\tau}{3} R_{ihjk}\Big)g_{lm}{{R_{ihjk}}} \\
&=\left(-R_{iabj}R_{kabh} R_{ihjk}  +R_{iabk}R_{jabh} R_{ihjk} +\frac{1}{2}R_{abih}R_{abjk} R_{ihjk} +\dfrac{1}{3}\tau R_{ihjk} R_{ihjk} \right)g_{lm}\\
&=\left(-2\mathring{R}+\hat{R}+\frac{1}{3}\tau||R||^2\right)g_{lm},
\end{aligned}
\end{equation*}
\begin{equation*}
\begin{aligned}
-&\Big(-T_{ijkl} +T_{ikjl}+\frac{1}{2}S_{iljk}+\frac{\tau}{3} R_{iljk}\Big)g_{hm}{{R_{ihjk}}}\\
&=\Big(R_{iabj}R_{kabl} R_{imjk} -R_{iabk}R_{jabl} R_{imjk}-\frac{1}{2}R_{abil}R_{abjk} R_{imjk}-\dfrac{\tau}{3} R_{iljk} R_{imjk}\Big)\\
&=2\mathring{R}_{lm}-\hat{R}_{lm}-\frac{\tau}{3}\tR_{lm},
\end{aligned}
\end{equation*}
\begin{equation*}
-\Big(-T_{ijmh}+T_{imjh}+\frac{1}{2}S_{ihjm}+\frac{\tau}{3} R_{ihjm}\Big)g_{lk}{{R_{ihjk}}}=2\mathring{R}_{lm}-\hat{R}_{lm}-\frac{\tau}{3}\tR_{lm},
\end{equation*}
\begin{align*}
&\Big(-T_{ijml} +T_{imjl}+\frac{1}{2}S_{iljm}+\frac{\tau}{3} R_{iljm}\Big)g_{hk}{{R_{ihjk}}}\\
&=\left(-R_{iabj}R_{mabl}+R_{iabm}R_{jabl}+\frac{1}{2}R_{abil}R_{abjm}+\frac{\tau}{3}R_{iljm}\right)(-\rho_{ij})\\
&=-\frac{\tau}{6}\left(-\frac{\tau^2}{36}g_{lm}+\tR_{lm}+\frac{1}{2}\tR_{lm}-\frac{\tau^2}{18}g_{lm}\right)\\
&=\frac{\tau^3}{72}g_{lm}-\frac{\tau}{4}\tR_{lm},
\end{align*}
\begin{equation*}
\Big(-T_{ikmh} +T_{imkh}+\frac{1}{2}S_{ihkm}+\frac{\tau}{3} R_{ihkm}\Big)g_{lj}{{ R_{ihjk}}}
=2\mathring{R}_{lm}-\hat{R}_{lm}-\frac{\tau}{3}\tR_{lm},
\end{equation*}
\begin{equation*}
-\Big(-T_{ikml}+T_{imkl}+\frac{1}{2}S_{ilkm}+\frac{\tau}{3} R_{ilkm}\Big)g_{hj}{{ R_{ihjk}}}
=\frac{\tau^3}{72}g_{lm}-\frac{\tau}{4}\tR_{lm},
\end{equation*}
\begin{equation*}
-\Big(-T_{hjml} +T_{hmjl}+\frac{1}{2}S_{hljm}+\frac{\tau}{3} R_{hljm}\Big)g_{ik}R_{ihjk}	=\frac{\tau^3}{72}g_{lm}-\frac{\tau}{4}\tR_{lm},
\end{equation*}
\begin{equation*}
\Big(-T_{hjkl} +T_{hkjl}+\frac{1}{2}S_{hljk}+\frac{\tau}{3} R_{hljk}\Big)g_{im}{{ R_{ihjk}}}
=2\mathring{R}_{lm}-\hat{R}_{lm}-\frac{\tau}{3}\tR_{lm},
\end{equation*}
\begin{equation*}
\Big(-T_{hkml} +T_{hmkl}+\frac{1}{2}S_{hlkm}+\frac{\tau}{3} R_{hlkm}\Big)g_{ij}{{ R_{ihjk}}} =\frac{\tau^3}{72}g_{lm}-\frac{\tau}{4}\tR_{lm}.
\end{equation*}
Thus, the third {{term}} of \eqref{eq:6-einstein(0,6)-tensor} transvecting with $R_{ihjk}$  becomes
\begin{equation*}
\left(-2\mathring{R}+\hat{R}+\frac{\tau}{3}||R||^2 + \frac{\tau^3}{18}\right)g_{lm}+8\mathring{R}_{lm}-4\hat{R}_{lm}-\frac{7}{3}\tau\tR_{lm}.
\end{equation*}
Making use of \eqref{eq:re}, the remaining terms in \eqref{eq:6-einstein(0,6)-tensor} transvecting with $R_{ihjk}$, we have
\begin{equation*}
\begin{aligned}	
&(A_{hjkmil}-A_{ljkmih}-A_{ijkmhl}+A_{ijmkhl}-A_{hjmkil}\\
&+A_{ljmkih}-A_{ikmjhl}+A_{hkmjil}-A_{lkmjih}){{ R_{ihjk}}}\\
=\;& R_{ahjk}R_{amil} R_{ihjk}-R_{aljk}R_{amih} R_{ihjk}-R_{aijk}R_{amhl} R_{ihjk}\\
&+R_{aijm}R_{akhl} R_{ihjk}-R_{ahjm}R_{akil} R_{ihjk}+R_{aljm}R_{akih} R_{ihjk}\\
&-R_{aikm}R_{ajhl} R_{ihjk}+R_{ahkm}R_{ajil} R_{ihjk}-R_{alkm}R_{ajih} R_{ihjk}\\
=&-4\check{R}_{lm}-2\hat{R}_{lm}+4\mathring{R}_{lm}.
\end{aligned}
\end{equation*}
{{Summing all terms, then we complete the proof of Theorem~B-\ref{eq:6-einstein}.

\smallskip

Now, we prove Theorem~B-\ref{eq:6-dim2stein}. Applying~\eqref{eq:supEin} to \eqref{eq:6-einstein(0,6)-tensor}, we have}}
	\begin{equation}\label{eq:6-2stein(0,6)-tensor}
		\begin{aligned}
			&-\frac{1}{8}\Big(||R||^2 -\frac{\tau^2}{3}\Big)\Big(g_{ij}g_{hk}g_{lm}-g_{ij}g_{hm}g_{lk}-g_{ik}g_{hj}g_{lm}+g_{ik}g_{hm}g_{lj}+g_{im}g_{hj}g_{lk}-g_{im}g_{hk}g_{lj}\Big)\\
			&+\Big(-T_{ijkh} + T_{ikjh} +\frac{1}{2}S_{ihjk}+\frac{\tau}{3} R_{ihjk}\Big)g_{lm}-\Big(-T_{ijkl} +T_{ikjl}+\frac{1}{2}S_{iljk}+\frac{\tau}{3} R_{iljk}\Big)g_{hm}\\
			&-\Big(-T_{ijmh}+T_{imjh}+\frac{1}{2}S_{ihjm}+\frac{\tau}{3} R_{ihjm}\Big)g_{lk}+\Big(-T_{ijml} +T_{imjl}+\frac{1}{2}S_{iljm}+\frac{\tau}{3} R_{iljm}\Big)g_{hk}\\
			&+\Big(-T_{ikmh} +T_{imkh}+\frac{1}{2}S_{ihkm}+\frac{\tau}{3} R_{ihkm}\Big)g_{lj}-\Big(-T_{ikml} +T_{imkl}+\frac{1}{2}S_{ilkm}+\frac{\tau}{3}R_{ilkm}\Big)g_{hj}\\
			&-\Big(-T_{hjml} +T_{hmjl}+\frac{1}{2}S_{hljm}+\frac{\tau}{3} R_{hljm}\Big)g_{ik}+\Big(-T_{hjkl} +T_{hkjl}+\frac{1}{2}S_{hljk}+\frac{\tau}{3} R_{hljk}\Big)g_{im}\\
			&+\Big(-T_{hkml} +T_{hmkl}+\frac{1}{2}S_{hlkm}+\frac{\tau}{3} R_{hlkm}\Big)g_{ij}+A_{hjkmil}-A_{ljkmih}-A_{ijkmhl}+A_{ijmkhl}\\
			&-A_{hjmkil}+A_{ljmkih}-A_{ikmjhl}+A_{hkmjil}-A_{lkmjih}=0.
		\end{aligned}
	\end{equation}
	{{We transvect {{each term}} of \eqref{eq:6-2stein(0,6)-tensor} with $R_{ihjk}$.}} Here, we give a few representative terms as follows:
		\begin{equation*}
			\begin{aligned}
				&-\frac{1}{8}\Big(||R||^2 -\frac{\tau^2}{3}\Big)\Big(g_{ij}g_{hk}g_{lm}-g_{ij}g_{hm}g_{lk}-g_{ik}g_{hj}g_{lm}+g_{ik}g_{hm}g_{lj}+g_{im}g_{hj}g_{lk}-g_{im}g_{hk}g_{lj}\Big){{R_{ihjk}}}\\&=\frac{1}{6}\tau\Big(||R||^2 -\frac{\tau^2}{3}\Big) g_{lm},
				\end{aligned}
		\end{equation*}
		\begin{equation*}
		    \Big(-T_{ijkh} + T_{ikjh} +\frac{1}{2}S_{ihjk}+\frac{\tau}{3} R_{ihjk}\Big)g_{lm}{{R_{ihjk}}} =\left(-2\mathring{R}+\hat{R}+\frac{8}{3}\tau||R||^2\right)g_{lm},
		\end{equation*}
		\begin{equation*}
		    -\Big(-T_{ijkl} +T_{ikjl}+\frac{1}{2}S_{iljk}+\frac{\tau}{3} R_{iljk}\Big)g_{hm}{{R_{ihjk}}}=2\mathring{R}_{lm}-\hat{R}_{lm}-\frac{1}{18}\tau||R||^2g_{lm},
		\end{equation*}
{{\begin{equation*}
\begin{aligned}	
&(A_{hjkmil}-A_{ljkmih}-A_{ijkmhl}+A_{ijmkhl}-A_{hjmkil}+A_{ljmkih}-A_{ikmjhl}+A_{hkmjil}-A_{lkmjih})R_{ihjk}\\
&=-4\check{R}_{lm}-2\hat{R}_{lm}+4\mathring{R}_{lm}\\
&=-\frac{\tau}{9}||R||^2 g_{lm}-2\hat{R}_{lm}+4\mathring{R}_{lm}.
\end{aligned}
\end{equation*}}}
Similarly, we can obtain remaining terms by transvecting with $R_{ihjk}$. Then, by rearranging all terms, we have Theorem~B-\ref{eq:6-dim2stein}.


\smallskip

Now we give an example of Theorem B.
\begin{ex}	
Let $M$ be a Riemannian product manifold of 3-dimensional Riemannian manifolds $M_1(k)$ and $M_2(k)$ of constant sectional curvature $k$ ($k \neq 0$). {{Then we can easily check that $M$ is an Einstein manifold. But $M$ can never be 2-stein: Let $X=(X_1 , X_2)$ is a tangent vector with $X_1$, $X_2$ its components tangent to $M_1 (k)$ and $M_2 (k)$, respectively. Then, $\Tr(R_X^2) = 2k^2 (||X_1 ||^4 + ||X_2 ||^4)$ which cannot be equal to $\mu_2(||X_1 ||^2 + ||X_2 ||^2)^2$.}} Let $\{e_{i}\}, \; i=1,2,\ldots,6$ be an orthonormal basis of $M$, where $\{e_{1},e_{2},e_{3}\}$ and $\{e_{4},e_{5},e_{6}\}$ are bases for $M_1(k)$ and $M_2(k)$, respectively. Then, we have
\begin{equation}\label{eq:ex_Ein}
R_{1221}=R_{1331}=R_{2332}=R_{4554}=R_{4664}=R_{5665}=k,
\end{equation}
and all the other components of $R$ vanish. From \eqref{eq:ex_Ein} we have
\begin{gather*}
		\tau = 12k,\quad ||R||^2 = 24 k^2,\quad \mathring{R} = -12k^3,\quad \hat{R} = -48k^3,\\
		\check{R}_{ij}=8k^3 \delta_{ij},\quad \hat{R}_{ij}=-8k^3 \delta_{ij},\quad \mathring{R}_{ij} = -2k^3 \delta_{ij},\quad \tR_{ij}=4k^2 \delta_{ij}.
		\end{gather*}
		Therefore, we find that the curvature identities of Theorem B hold on $M$. Here we note that $M$ is a super-Einstein manifold.
\end{ex}


\section{Appendix}
In this appendix we give the proof of Lemma \ref{thm:6-einstien(0,6)-tensor}.

{{By replacing}} the curvature tensor $R$ by the Weyl curvature tensor $W$ in (\ref{eq:Pa_n}), we can also obtain the curvature identity of $W$. In the case of $m=6$ and $r=2$, making use of the fact that the Weyl curvature tensor is traceless, we have the following curvature identity.
\begin{prop}The Weyl curvature tensor $W$ of any 6-dimensional Riemannian manifold satisfies the following identity{\rm :}
\begin{equation}\label{eq:weyl_iden}
    \begin{aligned}
    0=||W||^2\big\{&g_{ij}g_{hk}g_{lm}+g_{ik}g_{hm}g_{lj}+g_{im}g_{hj}g_{lk}-g_{ik}g_{hj}g_{lm}-g_{ij}g_{hm}g_{lk}-g_{im}g_{hk}g_{lj}\big\}\\
    -4\big\{&W_{iabc}W_{jabc}g_{hk}g_{lm}+W_{iabc}W_{kabc}g_{hm}g_{lj}+W_{iabc}W_{mabc}g_{hj}g_{lk}\\
    &-W_{iabc}W_{kabc}g_{hj}g_{lm}-W_{iabc}W_{jabc}g_{hm}g_{lk}-W_{iabc}W_{mabc}g_{hk}g_{lj}\\
    &-W_{habc}W_{jabc}g_{ik}g_{lm}-W_{habc}W_{kabc}g_{im}g_{lj}-W_{habc}W_{mabc}g_{ij}g_{lk}\\
    &+W_{habc}W_{kabc}g_{ij}g_{lm}+W_{habc}W_{jabc}g_{im}g_{lk}+W_{habc}W_{mabc}g_{ik}g_{lj}\\
    &+W_{labc}W_{jabc}g_{hk}g_{im}+W_{labc}W_{kabc}g_{hm}g_{ij}+W_{labc}W_{mabc}g_{hj}g_{ik}\\
    &-W_{labc}W_{kabc}g_{hj}g_{im}-W_{labc}W_{jabc}g_{hm}g_{ik}-W_{labc}W_{mabc}g_{hk}g_{ij}\big\}\\
    -8\big\{&(W_{iabj}W_{kabh}-W_{iabk}W_{jabh})g_{lm}-(W_{iabj}W_{mabh}-W_{iabm}W_{jabh})g_{lk}\\
    &+(W_{iabk}W_{mabh}-W_{iabm}W_{kabh})g_{lj}-(W_{iabj}W_{kabl}-W_{iabk}W_{jabl})g_{hm}\\
     &+(W_{iabj}W_{mabl}-W_{iabm}W_{jabl})g_{hk}-(W_{iabk}W_{mabl}-W_{iabm}W_{kabl})g_{hj}\\
    &+(W_{habj}W_{kabl}-W_{habk}W_{jabl})g_{im}-(W_{habj}W_{mabl}-W_{habm}W_{jabl})g_{ik}\\
    &+(W_{habk}W_{mabl}-W_{habm}W_{kabl})g_{ij}\big\}
    \end{aligned}
\end{equation}
\begin{align*}
    \quad\qquad\qquad\;+4\big\{&W_{abih}W_{abjk}g_{lm}-W_{abih}W_{abjm}g_{lk}+W_{abih}W_{abkm}g_{lj}-W_{abil}W_{abjk}g_{hm}\\
    &+W_{abil}W_{abjm}g_{hk}-W_{abil}W_{abkm}g_{hj}+W_{abhl}W_{abjk}g_{im}-W_{abhl}W_{abjm}g_{ik}\\
    &+W_{abhl}W_{abkm}g_{ij}\big\}\\
    \quad\qquad\qquad\;+8\big\{&W_{ahjk}W_{amil}-W_{aljk}W_{amih}-W_{ajhl}W_{aikm}-W_{aijk}W_{amhl}+W_{ajih}W_{almk}\\
    &+W_{ajil}W_{ahkm}-W_{ahjm}W_{akil}+W_{aljm}W_{akih}+W_{aijm}W_{akhl}\big\}.
\end{align*}
\end{prop}
Since the Weyl tensor W on a 6-dimensional Riemannian manifold is given by
\begin{equation*}
    W_{abcd} = R_{abcd} - \frac{1}{4}(\rho_{ad}g_{bc}+\rho_{bc}g_{ad} - \rho_{ac}g_{bd} -\rho_{bd}g_{ac}) + \frac{\tau}{20}(g_{ad}g_{bc} - g_{ac}g_{bd}),
\end{equation*}
we substitute the above Weyl tensor into \eqref{eq:weyl_iden} and  use the Einstein condition $\rho_{ij} = \frac{\tau}{6}g_{ij}$. Then, we obtain the following 34 terms.

{{\begin{enumerate}[label={(\arabic*)},leftmargin=0cm,align=left]
    \item $\big(||R||^2 - 4||\rho||^2 + \tau^2 \big)g_{ij}g_{hk}g_{lm}=\left(||R||^2 +\dfrac{\tau^2}{3} \right)g_{ij}g_{hk}g_{lm}$
    \item $-\big(||R||^2 - 4||\rho||^2 + \tau^2 \big)g_{ij}g_{hm}g_{lk}=-\left(||R||^2 +\dfrac{\tau^2}{3} \right)g_{ij}g_{hm}g_{lk}$
    \item $-\big(||R||^2 - 4||\rho||^2 + \tau^2 \big)g_{ik}g_{hj}g_{lm}=-\left(||R||^2 +\dfrac{\tau^2}{3} \right)g_{ik}g_{hj}g_{lm}$
    \item $\big(||R||^2 - 4||\rho||^2 + \tau^2 \big)g_{ik}g_{hm}g_{lj}=\left(||R||^2 +\dfrac{\tau^2}{3} \right)g_{ik}g_{hm}g_{lj}$
    \item $\big(||R||^2 - 4||\rho||^2 + \tau^2 \big)g_{im}g_{hj}g_{lk}=\left(||R||^2 +\dfrac{\tau^2}{3} \right)g_{im}g_{hj}g_{lk}$
    \item $-\big(||R||^2 - 4||\rho||^2 + \tau^2 \big)g_{im}g_{hk}g_{lj}=-\left(||R||^2 +\dfrac{\tau^2}{3} \right)g_{im}g_{hk}g_{lj}$
\smallskip
    \item $\big(-4R_{iabc}R_{jabc} + 8\rho_{ia}\rho_{ja} +8R_{iabj}\rho_{ab}-4\tau\rho_{ij}\big)g_{hk}g_{lm}=\left(-4R_{iabc}R_{jabc} -\dfrac{2}{9}\tau^2 g_{ij}\right)g_{hk}g_{lm}$
    \item $-\big(-4R_{iabc}R_{jabc} + 8\rho_{ia}\rho_{ja} +8R_{iabj}\rho_{ab}-4\tau\rho_{ij}\big)g_{hm}g_{lk}=-\left(-4R_{iabc}R_{jabc} -\dfrac{2}{9}\tau^2 g_{ij}\right)g_{hm}g_{lk}$
    \item $-\big(-4R_{iabc}R_{kabc} + 8\rho_{ia}\rho_{ka} +8R_{iabk}\rho_{ab}-4\tau\rho_{ik}\big)g_{hj}g_{lm}=-\left(-4R_{iabc}R_{kabc} -\dfrac{2}{9}\tau^2 g_{ik}\right)g_{hj}g_{lm}$
    \item $\big(-4R_{iabc}R_{kabc} + 8\rho_{ia}\rho_{ka} +8R_{iabk}\rho_{ab}-4\tau\rho_{ik}\big)g_{hm}g_{lj}=\left(-4R_{iabc}R_{kabc} -\dfrac{2}{9}\tau^2 g_{ik}\right)g_{hm}g_{lj}$
    \item $\big(-4R_{iabc}R_{mabc} + 8\rho_{ia}\rho_{ma} +8R_{iabm}\rho_{ab}-4\tau\rho_{im}\big)g_{hj}g_{lk}=\left(-4R_{iabc}R_{mabc} -\dfrac{2}{9}\tau^2 g_{im}\right)g_{hj}g_{lk}$
    \item $-\big(-4R_{iabc}R_{mabc} + 8\rho_{ia}\rho_{ma} +8R_{iabm}\rho_{ab}-4\tau\rho_{im}\big)g_{hk}g_{lj}=-\left(-4R_{iabc}R_{mabc} -\dfrac{2}{9}\tau^2 g_{im}\right)g_{hk}g_{lj}$
    \item $-\big(-4R_{habc}R_{jabc} + 8\rho_{ha}\rho_{ja} +8R_{habj}\rho_{ab}-4\tau\rho_{hj}\big)g_{ik}g_{lm}=-\left(-4R_{habc}R_{jabc} -\dfrac{2}{9}\tau^2 g_{hj}\right)g_{ik}g_{lm}$
    \item $\big(-4R_{habc}R_{jabc} + 8\rho_{ha}\rho_{ja} +8R_{habj}\rho_{ab}-4\tau\rho_{hj}\big)g_{im}g_{lk}=\left(-4R_{habc}R_{jabc} -\dfrac{2}{9}\tau^2 g_{hj}\right)g_{im}g_{lk}$
    \item $\big(-4R_{habc}R_{kabc} + 8\rho_{ha}\rho_{ka} +8R_{habk}\rho_{ab}-4\tau\rho_{hk}\big)g_{ij}g_{lm}=\left(-4R_{habc}R_{kabc} -\dfrac{2}{9}\tau^2 g_{hk}\right)g_{ij}g_{lm}$
    \item $-\big(-4R_{habc}R_{kabc} + 8\rho_{ha}\rho_{ka} +8R_{habk}\rho_{ab}-4\tau\rho_{hk}\big)g_{im}g_{lj}=-\left(-4R_{habc}R_{kabc} -\dfrac{2}{9}\tau^2 g_{hk}\right)g_{im}g_{lj}$
    \item $-\big(-4R_{habc}R_{mabc} + 8\rho_{ha}\rho_{ma} +8R_{habm}\rho_{ab}-4\tau\rho_{hm}\big)g_{ij}g_{lk}=-\left(-4R_{habc}R_{mabc} -\dfrac{2}{9}\tau^2 g_{hm}\right)g_{ij}g_{lk}$
    \item $\big(-4R_{habc}R_{mabc} + 8\rho_{ha}\rho_{ma} +8R_{habm}\rho_{ab}-4\tau\rho_{hm}\big)g_{ik}g_{lj}
        =\left(-4R_{habc}R_{mabc} -\dfrac{2}{9}\tau^2 g_{hm}\right)g_{ik}g_{lj}$
    \item $\big(-4R_{labc}R_{jabc} + 8\rho_{la}\rho_{ja} +8R_{labj}\rho_{ab}-4\tau\rho_{lj}\big)g_{ik}g_{hm}
        =\left(-4R_{labc}R_{jabc} -\dfrac{2}{9}\tau^2 g_{lj}\right)g_{ik}g_{hm}$
    \item $-\big(-4R_{labc}R_{jabc} + 8\rho_{la}\rho_{ja} +8R_{labj}\rho_{ab}-4\tau\rho_{lj}\big)g_{im}g_{hk}
        =-\left(-4R_{labc}R_{jabc} -\dfrac{2}{9}\tau^2 g_{lj}\right)g_{im}g_{hk}$
    \item $-\big(-4R_{labc}R_{kabc} + 8\rho_{la}\rho_{ka} +8R_{labk}\rho_{ab}-4\tau\rho_{lk}\big)g_{ij}g_{hm}
        =-\left(-4R_{labc}R_{kabc} -\dfrac{2}{9}\tau^2 g_{lk}\right)g_{ij}g_{hm}$
    \item $\big(-4R_{labc}R_{kabc} + 8\rho_{la}\rho_{ka} +8R_{labk}\rho_{ab}-4\tau\rho_{lk}\big)g_{im}g_{hj}
        =\left(-4R_{labc}R_{kabc} -\dfrac{2}{9}\tau^2 g_{lk}\right)g_{im}g_{hj}$
    \item $\big(-4R_{labc}R_{mabc} + 8\rho_{la}\rho_{ma} +8R_{labm}\rho_{ab}-4\tau\rho_{lm}\big)g_{ij}g_{hk}
        =\left(-4R_{labc}R_{mabc} -\dfrac{2}{9}\tau^2 g_{lm}\right)g_{ij}g_{hk}$
    \item $-\big(-4R_{labc}R_{mabc} + 8\rho_{la}\rho_{ma} +8R_{labm}\rho_{ab}-4\tau\rho_{lm}\big)g_{ik}g_{hj}
        =-\left(-4R_{labc}R_{mabc} -\dfrac{2}{9}\tau^2 g_{lm}\right)g_{ik}g_{hj}$
\smallskip
    \item \begin{align*}
        &\big(-8R_{iabj}R_{kabh} +8R_{iabk}R_{jabh} +4R_{abih}R_{abjk} -8R_{aijk}\rho_{ah} +8R_{akih}\rho_{aj} -8R_{ajih}\rho_{ak}\\
        &\quad+8R_{ahjk}\rho_{ai} +8\rho_{ij}\rho_{hk} -8\rho_{ik}\rho_{hj} -4\tau R_{ihjk}\big)g_{lm}\\
        &=\Big(-8R_{iabj}R_{kabh} +8R_{iabk}R_{jabh} +4R_{abih}R_{abjk} +\dfrac{4}{3}\tau R_{ihjk} +\dfrac{2}{9}\tau^2 g_{ij}g_{hk} - \dfrac{2}{9}\tau^2 g_{ik}g_{hj} \Big)g_{lm}
    \end{align*}
    \item \begin{align*}
        &-\big(-8R_{iabj}R_{kabl} +8R_{iabk}R_{jabl} +4R_{abil}R_{abjk} -8R_{aijk}\rho_{al} +8R_{akil}\rho_{aj} -8R_{ajil}\rho_{ak}\\
        &\quad+8R_{aljk}\rho_{ai} +8\rho_{ij}\rho_{lk} -8\rho_{ik}\rho_{lj} -4\tau R_{iljk}\big)g_{hm}\\
        &=-\Big(-8R_{iabj}R_{kabl} +8R_{iabk}R_{jabl} +4R_{abil}R_{abjk} +\dfrac{4}{3}\tau R_{iljk} +\dfrac{2}{9}\tau^2 g_{ij}g_{lk} - \dfrac{2}{9}\tau^2 g_{ik}g_{lj}\Big)g_{hm}
    \end{align*}
    \item \begin{align*}
        &-\big(-8R_{iabj}R_{mabh} +8R_{iabm}R_{jabh} +4R_{abih}R_{abjm} -8R_{aijm}\rho_{ah} +8R_{amih}\rho_{aj} -8R_{ajih}\rho_{am}\\
        &\quad+8R_{ahjm}\rho_{ai} +8\rho_{ij}\rho_{hm} -8\rho_{im}\rho_{hj} -4\tau R_{ihjm}\big)g_{lk}\\
        &=-\Big(-8R_{iabj}R_{mabh} +8R_{iabm}R_{jabh} +4R_{abih}R_{abjm} +\dfrac{4}{3}\tau R_{ihjm} +\dfrac{2}{9}\tau^2 g_{ij}g_{hm} - \dfrac{2}{9}\tau^2 g_{im}g_{hj}\Big)g_{lk}
    \end{align*}
    \item \begin{align*}
        &\big(-8R_{iabj}R_{mabl} +8R_{iabm}R_{jabl} +4R_{abil}R_{abjm} -8R_{aijm}\rho_{al} +8R_{amil}\rho_{aj} -8R_{ajil}\rho_{am}\\
        &\quad+8R_{aljm}\rho_{ai} +8\rho_{ij}\rho_{lm} -8\rho_{im}\rho_{lj} -4\tau R_{iljm}\big)g_{hk}\\
    &=\Big(-8R_{iabj}R_{mabl} +8R_{iabm}R_{jabl} +4R_{abil}R_{abjm} +\dfrac{4}{3}\tau R_{iljm} +\dfrac{2}{9}\tau^2 g_{ij}g_{lm} - \dfrac{2}{9}\tau^2 g_{im}g_{lj}\Big)g_{hk}
    \end{align*}
    \item \begin{align*}
        &\big(-8R_{iabk}R_{mabh} +8R_{iabm}R_{kabh} +4R_{abih}R_{abkm} -8R_{aikm}\rho_{ah} +8R_{amih}\rho_{ak} -8R_{akih}\rho_{am} \\
        &\quad+8R_{ahkm}\rho_{ai} +8\rho_{ik}\rho_{hm} -8\rho_{im}\rho_{hk} -4\tau R_{ihkm}\big)g_{lj}\\
    &=\Big(-8R_{iabk}R_{mabh} +8R_{iabm}R_{kabh} +4R_{abih}R_{abkm} +\dfrac{4}{3}\tau R_{ihkm} +\dfrac{2}{9}\tau^2 g_{ik}g_{hm} - \dfrac{2}{9}\tau^2 g_{im}g_{hk}\Big)g_{lj}
    \end{align*}
    \item \begin{align*}
        &-\big(-8R_{iabk}R_{mabl} +8R_{iabm}R_{kabl} +4R_{abil}R_{abkm} -8R_{aikm}\rho_{al} +8R_{amil}\rho_{ak} -8R_{akil}\rho_{am} \\
        &\quad+8R_{alkm}\rho_{ai} +8\rho_{ik}\rho_{lm} -8\rho_{im}\rho_{lk} -4\tau R_{ilkm}\big)g_{hj}\\
    &=-\Big(-8R_{iabk}R_{mabl} +8R_{iabm}R_{kabl} +4R_{abil}R_{abkm} +\dfrac{4}{3}\tau R_{ilkm} +\dfrac{2}{9}\tau^2 g_{ik}g_{lm} - \dfrac{2}{9}\tau^2 g_{im}g_{lk}\Big)g_{hj}
    \end{align*}
    \item \begin{align*}
        &-\big(-8R_{habj}R_{mabl} +8R_{habm}R_{jabl} +4R_{abhl}R_{abjm} -8R_{ahjm}\rho_{al} +8R_{amhl}\rho_{aj} -8R_{ajhl}\rho_{am}\\
        &\quad+8R_{aljm}\rho_{ah} +8\rho_{hj}\rho_{lm} -8\rho_{hm}\rho_{lj} -4\tau R_{hljm}\big)g_{ik}\\
    &=-\Big(-8R_{habj}R_{mabl} +8R_{habm}R_{jabl} +4R_{abhl}R_{abjm} +\dfrac{4}{3}\tau R_{hljm} +\dfrac{2}{9}\tau^2 g_{hj}g_{lm} - \dfrac{2}{9}\tau^2 g_{hm}g_{lj}\Big)g_{ik}
    \end{align*}
    \item \begin{align*}
        &\big(-8R_{habj}R_{kabl} +8R_{habk}R_{jabl} +4R_{abhl}R_{abjk} -8R_{ahjk}\rho_{al} +8R_{akhl}\rho_{aj} -8R_{ajhl}\rho_{ak} \\
        &\quad+8R_{aljk}\rho_{ah} +8\rho_{hj}\rho_{lk} -8\rho_{hk}\rho_{lj} -4\tau R_{hljk}\big)g_{im}\\
    &=\Big(-8R_{habj}R_{kabl} +8R_{habk}R_{jabl} +4R_{abhl}R_{abjk} +\dfrac{4}{3}\tau R_{hljk} +\dfrac{2}{9}\tau^2 g_{hj}g_{lk} - \dfrac{2}{9}\tau^2 g_{hk}g_{lj}\Big)g_{im}
    \end{align*}
    \item \begin{align*}
        &\big(-8R_{habk}R_{mabl} +8R_{habm}R_{kabl} +4R_{abhl}R_{abkm} -8R_{ahkm}\rho_{al} +8R_{amhl}\rho_{ak} -8R_{akhl}\rho_{am} \\
        &\quad+8R_{alkm}\rho_{ah} +8\rho_{hk}\rho_{lm} -8\rho_{hm}\rho_{lk} -4\tau R_{hlkm}\big)g_{ij}\\
    &=\Big(-8R_{habk}R_{mabl} +8R_{habm}R_{kabl} +4R_{abhl}R_{abkm} +\dfrac{4}{3}\tau R_{hlkm} +\dfrac{2}{9}\tau^2 g_{hk}g_{lm} - \dfrac{2}{9}\tau^2 g_{hm}g_{lk}\Big)g_{ij}
    \end{align*}
\smallskip
    \item \begin{align*}
        &8\Big(R_{ahjk}R_{amil}-R_{aljk}R_{amih}-R_{aijk}R_{amhl}+R_{aijm}R_{akhl}-R_{ahjm}R_{akil}+R_{aljm}R_{akih}\\
        &\;-R_{aikm}R_{ajhl}+8R_{ahkm}R_{ajil}-R_{alkm}R_{ajih} -R_{ilkm}\rho_{hj}+R_{iljm}\rho_{hk}+R_{hljk}\rho_{im}\\
        &\;+R_{ihjk}\rho_{lm}+R_{ihkm}\rho_{lj}-R_{ihjm}\rho_{lk}-R_{iljk}\rho_{hm}+R_{hlkm}\rho_{ij}-R_{hljm}\rho_{ik}\Big)\\
        &=8\big(R_{ahjk}R_{amil}-R_{aljk}R_{amih}-R_{aijk}R_{amhl}+R_{aijm}R_{akhl}-R_{ahjm}R_{akil}+R_{aljm}R_{akih}\\
        &\;-R_{aikm}R_{ajhl}+R_{ahkm}R_{ajil}-R_{alkm}R_{ajih}\Big)-\dfrac{4}{3}\tau \Big(R_{ilkm}g_{hj}- R_{iljm}g_{hk}- R_{hljk}g_{im}\\
        &\;- R_{ihjk}g_{lm}- R_{ihkm}g_{lj}+ R_{ihjm}g_{lk}+ R_{iljk}g_{hm}- R_{hlkm}g_{ij}+ R_{hljm}g_{ik}\Big).
    \end{align*}
\end{enumerate}}}
By rearranging all terms, we complete the proof of Lemma~\ref{thm:6-einstien(0,6)-tensor}.
\vskip 0.3cm
\section*{{Acknowledgements}}
The authors thank Prof. Nikolayevsky for several useful discussions.

\end{document}